\documentclass[12pt]{amsart}

\usepackage{amsmath}
\usepackage{amssymb}
\usepackage{amsthm}
\usepackage{latexsym}
\usepackage{graphicx}
\usepackage{cite}

\theoremstyle{plain}

\newtheorem{theo}{Theorem}

\newtheorem{prop}{Proposition}

\newtheorem*{fprop}{Finiteness property}
\newtheorem{conj}{Conjecture}

\theoremstyle{definition}

\theoremstyle{remark}

\begin{document}

\thispagestyle{empty}

\date{\today}
%

\thanks{The research reported here
was performed while the authors were at MIT, Cambridge. It was
partially supported by the ``Communaut\'e francaise de Belgique -
Actions de Recherche Concert\'ees", by the EU HYCON Network of
Excellence (contract number FP6-IST-511368), and by the Belgian
Programme on Interuniversity Attraction Poles initiated by the
Belgian Federal Science Policy Office
The scientific responsibility rests with its authors. Rapha\"el
Jungers is a FNRS fellow (Belgian Fund for Scientific Research).}

\title{On the finiteness property for rational matrices}

\author{
Rapha\"el Jungers}
\address{Division of Applied Mathematics,
Universit\'e catholique de Louvain, 4 avenue Georges Lemaitre,
B-1348 Louvain-la-Neuve, Belgium} \email{blondel@inma.ucl.ac.be}

\author{
Vincent D. Blondel}
\address{Division of Applied Mathematics,
Universit\'e catholique de Louvain, 4 avenue Georges Lemaitre,
B-1348 Louvain-la-Neuve, Belgium} \email{jungers@inma.ucl.ac.be}

%
%


\maketitle
Abstract. {We analyze the periodicity of optimal long products of matrices.  A set of matrices is said to have the finiteness property if the maximal rate of growth of long products of matrices taken from the set can be obtained by a periodic product.  It was conjectured a decade ago that all finite sets of real matrices have the finiteness property. This conjecture, known as the ``finiteness conjecture", is now known to be false but no explicit counterexample to the conjecture is available and in particular it is unclear if a counterexample is possible whose matrices have rational or binary entries.    In this paper, we prove that finite sets of nonnegative rational matrices have the finiteness property if and only if \emph{pairs} of \emph{binary} matrices do. We also show that all {pairs} of $2 \times 2$ binary matrices have the finiteness property. These results have direct implications for the stability problem for sets of matrices. Stability is algorithmically decidable for sets of matrices that have the finiteness property and so it follows from our results that if all pairs of binary matrices have the finiteness property then stability is decidable for sets of nonnegative rational matrices. This would be in sharp contrast with the fact that the related problem of boundedness is known to be undecidable for sets of nonnegative rational matrices.  }\\

\begin{section}{Introduction}
The joint spectral radius of a set of matrices characterizes the maximum rate of growth that can be obtained by forming long products of matrices.  Let $\Sigma \subset \mathbb{R}^{n \times n}$ be a finite set of matrices. The \emph{joint spectral radius} of  $\Sigma$ is defined by:
\begin{eqnarray}\label{def-jsr}\rho(\Sigma) &= &\limsup_{t\to \infty}
\max\{\|A\|^{1/t}: A  \in \, \Sigma^{t}\}\end{eqnarray}
 where $\Sigma^{t}$ is the set of products of length $t$ of matrices from $\Sigma$, i.e., $\Sigma^t=\{A_{1}\dots A_{t}:A_i\in\Sigma\}$.  It is easy to verify that the quantity $\rho(\Sigma)$ does not depend on the chosen matrix norm. It has been proved in \cite{berger-wang} that the following equality holds for finite (or bounded) sets $\Sigma$:
 \begin{eqnarray}\label{def-jsr2}\rho(\Sigma) &= &\limsup_{t\to \infty}
\nonumber\max\{\rho(A)^{1/t}: A  \in \, \Sigma^t\}
\end{eqnarray}
($\rho$ is used here to denote the usual spectral radius). It is also known that the following inequalities hold for all $t$:
\begin{eqnarray}\label{a1}\max{\{\rho(A)^{1/t}: A\in \Sigma^t\}}\leq\rho(\Sigma)\leq \max{\{||A||^{1/t}: A\in \Sigma^t\}}.\end{eqnarray}
  These inequalities provide a straightforward  way to approximate the joint spectral radius to any desired accuracy: evaluate the upper and lower bounds for products of increasing length $t$, until $\rho$ is squeezed in a sufficiently small interval and the desired accuracy is reached. Unfortunately, this method, and in fact any other general method for computing or approximating the joint spectral radius, is  bound to be inefficient.  Indeed, it has been proved that, unless $P=NP$, there is no algorithm that computes or even approximates with a priori guaranteed accuracy the joint spectral radius of a set of matrices in polynomial time \cite{tsitsiklis97lyapunov}.  And this is true even if the matrices have binary entries.

 For some sets $\Sigma$, the right hand side inequality in (\ref{a1}) is strict for all $t$.  This is the case for example for the set consisting of just one matrix $$\begin{pmatrix} 1&1\\0&1 \end{pmatrix}.$$
Thus, there is no hope to reach the exact value of the joint spectral radius by simply evaluating the right hand side in (\ref{a1}). On the other hand, since $\rho(A^k)=\rho^k(A)$ the left hand side always provides the exact value when the set  $\Sigma$ consists of only one matrix and one can  thus hope to reach the exact value of the joint spectral radius by  evaluating the maximal spectral radii of products of increasing length. If for some $t$ and  $A\in \Sigma^t$ we have $\rho(A)^{1/t}=\rho(\Sigma)$, then the value of the joint spectral radius is reached. Sets of matrices for which such a product is possible are said to have the finiteness property.
\begin{fprop}
A set $\Sigma$ of matrices is said to have the \emph{finiteness property} if there exists some product $A=A_1 \ldots A_t$ with $A_i \in \Sigma$  for which
$ \rho(\Sigma)=\rho^{1/t}(A).$
\end{fprop}

One of the interests of the finiteness property arises from its connection with the stability question for a set of matrices. This  problem is of practical interest in a number of application contexts. A set of matrices $\Sigma$ is \emph{stable} if all long products of matrices taken from the set converge to zero. There are no known algorithms for deciding stability of a set of matrices and it is unknown if this problem is algorithmically decidable. One can verify that stability of the set $\Sigma$ is equivalent to the condition $\rho(\Sigma)<1$ and we may therefore hope to decide stability as follows: for increasing values of $t$ evaluate $\underline {\rho}_t=\max\{\rho(A)^{1/t}: A\in \Sigma^t\}$ and $\overline {\rho}_t=\max\{||A||^{1/t}: A\in \Sigma^t\}$. From (\ref{a1}) we know that $\underline {\rho}_t\leq \rho \leq \overline {\rho}_t$ and so, as soon as a $t$ is reached for which $\overline {\rho}_t<1$ we stop and declare the set stable, and as soon as a $t$ is reached for which   $\underline {\rho}_t \geq 1$ we stop and declare the set unstable. This algorithm will always stop unless $\rho=1$ and $\underline {\rho}_t < 1$ for all $t$. But this last situation never occurs for sets of matrices that satisfy the finiteness property and so we conclude:

\begin{prop}
Stability is algorithmically decidable for sets of matrices that have the finiteness property.
\end{prop}

It was first conjectured in 1995 by  Lagarias and Wang that all sets of  real matrices have the finiteness property \cite{lagarias-finiteness}. This conjecture, known as the finiteness conjecture, has attracted intense attention in recent years and has  been proved to be false by several authors \cite{cfkoz, cfbousch, blondel-elementary}. So far all proofs provided are nonconstructive and all sets of matrices whose joint spectral radius is known exactly satisfy the finiteness property. The finiteness property is also known to hold in a number of particular cases including the case were the matrices are symmetric, or if the Lie algebra associated with the set of matrices is solvable \cite[Corollary 6.19]{theys_thesis}: in this case the joint spectral radius is simply equal to the maximum of the spectral radii of the matrices.  This follows directly from a result of Liberzon et al. \cite{lib} after a conjecture of Gurvits \cite{Gu1}.

The definition of the finiteness property leads to a number of natural questions: When does the finiteness property holds? Is is decidable to determine if a given set of matrices satisfies the finiteness property?  Do  matrices with  rational entries satisfy the finiteness property? Do  matrices with binary entries satisfy the finiteness property?  In the first theorem in this paper we prove a connection between rational and binary matrices:

\begin{theo}
The finiteness property holds for all sets of nonnegative rational  matrices if and only if it holds for  all \emph{pairs} of \emph{binary} matrices.
\end{theo}

The case of binary matrices  appears to be important in a number applications \cite{mos-bounds, MoisionOrlitskySiegel01, moision-joint-constraints,Crespi05, protasov3}.  These applications have led to a number of joint spectral radius computations for binary matrices  \cite{mos-bounds, MoisionOrlitskySiegel01, master_thesis}.  The results obtained so far seem to indicate that for binary matrices there is always an optimal infinite periodic product.  Moreover, when the matrices have binary entries  they can be interpreted as adjacency matrices of  graphs on an identical set of nodes and in this context it seems  natural to expect optimality to be obtained for periodic products. Motivated by these observations, the following conjecture is phrased in \cite{BlondelJungersProtasov06}:

\begin{conj}\label{conj-B}
Pairs of binary matrices have the finiteness property.
\end{conj}

Of course if this conjecture is correct then  nonnegative rational matrices also satisfy the finiteness property and this in turn implies that stability, that is, the question $\rho <1$,  is decidable for sets of matrices with nonnegative rational entries. From a decidability perspective this last result would be somewhat surprising since it is known that the closely related question $\rho \leq 1$ is not algorithmically decidable for such sets of matrices \cite{blondel03undecidable, blondel00boundedness}.

Motivated by the relation between binary and rational matrices, we prove in our second theorem that pairs of $2 \times 2$ binary matrices satisfy the finiteness property. We have not been able to find a unique argument for all possible pairs  and we therefore proceed by enumerating a number of  cases and by providing separate  proofs for each of them. This somewhat unsatisfactory proof is nevertheless encouraging in that it could possibly be representative of the difficulties arising for pairs of binary matrices of arbitrary dimension.  In particular, some of the techniques we use can be applied to matrices of arbitrary dimension.

Let us finally notice that in all the numerical computations that we have performed on binary matrices not only the finiteness property always seemed to hold but the  period length of optimal products was always very short.  The computation of the joint spectral radius is known to be NP-hard for binary matrices but this does not exclude the possibility of a bound on the period length that is  linear in the dimensions of the matrices.  In the case of matrices characterizing the capacity of codes avoiding forbidden difference patterns, the length of the period is even suspected to be sublinear (see Conjecture 1 in \cite{master_thesis}).

\end{section}

\begin{section}{Finiteness property for rational and binary matrices}

In this section, we prove that the finiteness property holds for nonnegative rational matrices if and only if it holds for  \emph{pairs} of \emph{binary} matrices. The proof proceeds in three steps.  First we reduce the nonnegative rational case to the nonnegative integer case, we then reduce this case to the binary case, and finally we show how to reduce the number of matrices to two.

\begin{prop}\label{lemQZ}
The finiteness property holds for sets of nonnegative rational matrices  if and only if it holds for sets of nonnegative integer matrices.
\end{prop}

\begin{proof}
From the definition of the joint spectral radius we have that for any $\alpha >0,\ \rho(\Sigma)= (1/\alpha) \rho(\alpha \Sigma)$.  Now, for any set $\Sigma$ of matrices with nonnegative rational entries , let us pick an $\alpha \neq 0 \in \mathbb{N}$ such that $\alpha \Sigma \subseteq \mathbb{N}^{ n \times  n}$. If there exists a positive integer $t$ and a matrix $A\in (\alpha\Sigma)^t$ such that $\rho (\alpha\Sigma)=\rho^{1/t} (A)$, then $\rho(\Sigma)=(1/\alpha)\rho^{1/t}(A)=\rho^{1/t}(A/\alpha^t)$, where $A/\alpha^t\in \Sigma^t$.
\end{proof}
We now turn to the reduction from the integer to the binary case.
\begin{theo}\label{lemZB}
The finiteness property holds for sets of nonnegative integer matrices  if and only if it holds for sets of binary matrices.
\end{theo}\begin{proof}
Consider a finite set of nonnegative integer matrices $\Sigma \subset \mathbb{N}^{n \times n}$.  We think of the matrices in $\Sigma$ as  adjacency matrices of weighted graphs on a set of $n$ nodes and we construct auxiliary graphs such that paths of weight $w$ in the original weighted graphs are replaced by $w$ paths of weight one in the auxiliary graphs.  For every matrix $A\in \Sigma \subset \mathbb{N}^{n \times n}$, we introduce a new matrix $\tilde{A} \in  \{0,1\}^{nm\times nm}$  as follows.  We define $m$ as the largest entry of the matrices in $\Sigma$. Then, for every node $v_i$ $(i=1, \ldots, n)$ in the original graphs, we introduce $m$ nodes $\tilde{v}_{i,1},\ldots, \tilde{v}_{i,m}$ in the auxiliary graphs. The auxiliary graphs have $nW$ nodes; we now define their edges. For all $A\in \Sigma$ and $A_{i,j}=k\neq 0$, we define $km$ edges in $\tilde{A}$ from nodes $\tilde{v}_{i,s}:1\leq s \leq k$ to the nodes $\tilde{v}_{j,t}: 1\leq t \leq m$. The following reasoning leads now to the claim:
\begin{itemize}
\item For any product $A\in \Sigma^t$, and any couple of indices $(i,j)$, the corresponding product $\tilde{A}\in \tilde{\Sigma}^t$ has the following property: $\forall s,A_{i,j}=\sum_r\tilde{A}_{\tilde{v}_{i,r},\tilde{v}_{j,s}}$.  This is easy to show by induction on the length of the product.
\item  $\forall t,\forall A\in \Sigma^{t},\ ||A||_{1}=||\tilde{A}||_1$, where $||\cdot||_1$ represents the maximum sum of the absolute values of all entries of any column in a matrix.
\item $\rho(\Sigma)=\rho(\tilde{\Sigma})$, and if $\rho(\tilde{\Sigma})=\rho^{1/T}(\tilde{A}):\tilde{A}\in \tilde \Sigma^T$, then $\rho(\Sigma)=\rho^{1/T}(A)$, where $A$ is the product in $\Sigma^T$ corresponding to $\tilde A$.\end{itemize}
\end{proof}

We finally consider the last reduction: we are given a set of matrices and we reformulate the finiteness property for this set into the  finiteness property for two particular matrices constructed from the set. The construction is such that the entries of the two matrices are exactly those of the original matrices, except for some entries that are  equal to zero or one.

More specifically, assume that we are given $m$ matrices $A_1, \ldots, A_m$ of dimension $n$. From these $m$ matrices we construct two matrices $\tilde{A}_0,\tilde{A}_1$ of dimension $(2m-1)n$. The matrices $\tilde{A}_0,\tilde{A}_1$ consist of $(2m-1) \times (2m-1)$ square blocks of dimension $n$ that are either equal to $0$, $I$ or to one of the matrices $A_i$. The explicit construction of these two matrices is best illustrated with a graph.

Consider the graph $G_0$ on a set of $2m-1$ nodes $s_i$ $(i=1, \ldots, 2m-1)$ and  whose edges are given by $(i, i+1)$ for $i=1, \ldots, 2m-2$.  We also consider a graph $G_1$ defined on the same set of nodes and whose edges of weight $a_i$ are given by $(m+i-1, i)$ for $i=1, \ldots, m$. These two graphs are represented on Figure \ref{fig-subspaces} for the case $m=5$. In such a graph a directed path that leaves the node $m$ returns there after $m$ steps and whenever it does so, the path passes exactly once through an edge of graph $G_1$.  Let us now describe how to construct the matrices $\tilde{A}_0,\tilde{A}_1$. The matrices are obtained by constructing the adjacency matrices of the graphs $G_0$ and $G_1$ and by replacing the entries $1$ and $0$ by the matrices $I$ and $0$ of dimension $n$, and the weight $a_i$ by the matrices $A_i$. For the case $m=5$ the  matrices $\tilde{A}_0,\tilde{A}_1$ are thus given by:

$$\tilde{A}_0 = \begin{pmatrix}
0&I&0&0&0&0&0&0&0\\
0&0&I&0&0&0&0&0&0\\
0&0&0&I&0&0&0&0&0\\
0&0&0&0&I&0&0&0&0\\
0&0&0&0&0&I&0&0&0\\
0&0&0&0&0&0&I&0&0\\
0&0&0&0&0&0&0&I&0\\
0&0&0&0&0&0&0&0&I\\
0&0&0&0&0&0&0&0&0
\end{pmatrix} \;
\tilde{A}_1 =
\begin{pmatrix}
0&0&0&0&0&0&0&0&0\\
0&0&0&0&0&0&0&0&0\\
0&0&0&0&0&0&0&0&0\\
0&0&0&0&0&0&0&0&0\\
A_1&0&0&0&0&0&0&0&0\\
0&A_2&0&0&0&0&0&0&0\\
0&0&A_3&0&0&0&0&0&0\\
0&0&0&A_4&0&0&0&0&0\\
0&0&0&0&A_5&0&0&0&0
\end{pmatrix}$$

\begin{figure}
\includegraphics[width=0.4\textwidth ]{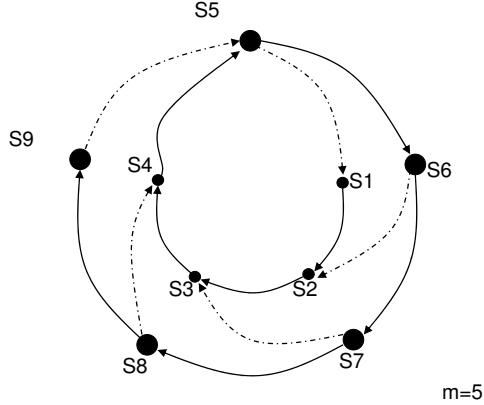} 
 \caption{Schematic representation of the macro transitions between subspaces.  The full edges represent transitions in $\tilde{A}_0$ and the dashed edges transitions in $\tilde{A}_1$.}
\label{fig-subspaces}
 \end{figure}

The two matrices so constructed inherit some of the properties of the graphs $G_0$ and $G_1$ and this allows us to prove the following theorem.

\begin{theo}\label{lem2matrices}
Consider a set of $m\geq 1$ matrices $\Sigma= \{A_1, \ldots, A_m: A_i \in \mathbb{R}^{n\times n}\}$ and define $\tilde \Sigma=\{\tilde{A}_0,\tilde{A}_1\}$ with the matrices $\tilde{A}_0$ and $\tilde{A}_1$ as defined above. Then $\rho(\tilde \Sigma)=\rho(\Sigma)^{1/m}$. Moreover, the finiteness property holds for $\Sigma$ if and only it holds for $\tilde \Sigma$.
\end{theo}

\begin{proof}

The crucial observation for the proof is the following. Consider a path in $G_0$ and $G_1$. Edges in $G_0$ and $G_1$ have outdegree at most equal to one and so if a sequence of graphs among $G_0$ and $G_1$ is given, there is only one path leaving $i$ that follows that particular sequence. This fact ensures that any block in any product of matrices in $\tilde{\Sigma}$ is a pure product of blocks of the matrices in $\tilde{\Sigma}$, and not a sum of such products. Moreover, any path leaving from $i$ and of length  $km$  either returns to $i$ after passing through $k$ edges of $G_1$, or ends at node $i+m$ after passing through $k-1$ edges of $G_1$, or ends at node $i+m$ (mod $2m$) after passing through $k+1$ edges of $G_1$. From this it follows that in a product of length $km$ of the matrices $\tilde{A}_0$ and $\tilde{A}_1$ there is exactly one nonzero block in every line of blocks, and this block is a product of length $k-1$, $k$, or $k+1$ of matrices from $\Sigma$.

We now show that $\rho(\tilde{\Sigma})\geq\rho(\Sigma)^{1/m}$ by proving that for any matrix $A\in \Sigma^{t}$, there is a matrix $\tilde{A} \in \tilde \Sigma^{tm}$ such that $||\tilde A||\geq ||A||$. Define  $\tilde{B}_i=\tilde{A}_{0}^{i-1}\tilde{A}_{1}\tilde{A}_{0}^{m-i}\in \tilde{\Sigma}^{m}$ for $i=1, \ldots, m$ so that the block in position $(m,m)$ in $\tilde{B}_i$ is simply equal to $A_i$. Consider now  some product of length $t$, $A=A_{i_1} \cdots A_{i_t} \in \Sigma^{t}$ and construct the corresponding matrix product $\tilde{A}=\tilde{B}_{i_1}\dots \tilde{B}_{i_t}\in \tilde{\Sigma}^{tm}$. The block in position $(m,m)$ in $\tilde A$ is equal to $A_{i_1}\dots A_{i_t}$ and so  $||\tilde A||\geq ||A||$ and $\rho(\tilde{\Sigma})\geq\rho(\Sigma)^{1/m}$.

Let us now show that $\rho(\tilde{\Sigma})\leq \rho(\Sigma)^{1/m}$. Consider therefore an arbitrary product $\tilde A \in \tilde \Sigma^{l}$ and decompose $\tilde A=\tilde C \tilde A'$ with $\tilde C$ a product of at most $m$ factors and $\tilde A' \in \Sigma^{km}$. By the observation above we know that there is at most one nonzero block in every line of blocks of $\tilde A'$, and this block is a product of length $k-1$, $k$, or $k+1$ of matrices from $\Sigma$. Therefore, if the norm is chosen to be the maximum line sum, we have $||\tilde A|| \leq  K_1 K_2 ||A||$ where $A$ is some product of length $k-1$ of matrices from $\Sigma$, $K_1$ is the maximal norm of a product of at most $m$ matrices in $\tilde \Sigma$, and $K_2$ is the maximal norm of a product of at most $2$ matrices in $\Sigma$. Using the last inequality, we arrive at $$||\tilde A||^{1/(k-1)} \leq  (K_1 K_2)^{1/(k-1)} ||A||^{1/(k-1)}.$$
The initial product $\tilde A$ is an arbitrary product of length $l=km+r$ and so by letting $k$ tend to infinity and using the definition of the joint spectral radius we conclude $\rho(\tilde{\Sigma})\leq \rho(\Sigma)^{1/m}$.

We have thus proved that $\rho(\tilde \Sigma)=\rho(\Sigma)^{1/m}$. It remains to prove the equivalence of the finiteness property.  If $\Sigma$ satisfies the finiteness property then $\rho(\Sigma)=\rho(A_1\dots A_T)^{1/t}$, then $\rho (\tilde \Sigma)=\rho(\Sigma)^{1/m}=\rho(\tilde B_1\dots \tilde B_T)^{1/(Tm)}$ and so $\tilde \Sigma$ does too. In the opposite direction, if the finiteness property holds for $\tilde \Sigma$, then we must have $\rho(\tilde \Sigma)=\rho(\tilde B_1\dots \tilde B_T)^{1/T}$, and then $\rho(\Sigma)=\rho(\tilde \Sigma )^{m}=\rho(A_1\dots A_T)^{1/T}$.

\end{proof}

Combining the results obtained so far in this section, we now state our main result. Before we do so, let us notice that the 

\begin{theo}
The finiteness property holds for all sets of nonnegative rational  matrices if and only if it holds for  all {pairs} of {binary} matrices.

The finiteness property holds for all sets of  rational  matrices if and only if it holds for  all {pairs} of  matrices with entries in $\{0,1,-1\}$.
\end{theo}

\begin{proof}
The proof for the nonnegative case is a direct consequence of Proposition \ref{lemQZ}, Theorem \ref{lemZB} and Theorem \ref{lem2matrices}.
For the  case of arbitrary rational entries, the results and proofs of Proposition \ref{lemQZ} and Theorem \ref{lem2matrices} may be applied as they are. For the construction in the proof of Theorem \ref{lemZB}, one just has to weight any edge with $-1$ whenever it represents a negative entry, and the reasoning in the proof still holds.
\end{proof}
Let us finally remark that for the purpose of reducing the finiteness property of rational matrices to pairs of binary matrices, we have proved here that for any set $\Sigma$ of $m$ matrices with integer (resp. nonnegative integer) entries, there is always a \emph{pair} of matrices $\tilde{\Sigma}$ with entries in $\{0,1,-1\}$ (resp. in $\{0,1\}$) such that $\rho(\Sigma)=\rho(\tilde\Sigma)^{m}$.  Loosely speaking pairs of binary matrices have the same combinatorial complexity.  
\end{section}

\begin{section}{The finiteness property for pairs of $2\times 2 $ binary matrices.}\label{Section-2X2}

In this section, we prove that the finiteness property holds for {pairs} of binary matrices of size $2 \times 2$.  Even if this result may at first sight appear anecdoctic, it has some relevance since it has been shown in the previous section that pairs of binary matrices are not restrictive.  Moreover, even for this $2 \times 2$ case, non-trivial behaviours occur.  For instance, the set of matrices $$\left \{\begin{pmatrix}1&1\\0&1 \end{pmatrix},\begin{pmatrix}0&1\\1&0 \end{pmatrix}  \right \}, $$ whose behaviour could at first sight seem very simple, happens to have a joint spectral radius equal to $((3+\sqrt{13})/2)^{1/4}$, and this value is only reached by products of length at least four.  Another interest of this section is to present techniques that may prove useful to establish the finiteness property for matrices of larger dimension.

There are $256$ ordered pairs of binary matrices.  Since we are only interested in unordered sets we can lower this number to $(2^4(2^4-1))/2=120$.  We first present a series of simple properties that allow us to handle most of these cases and we then give a complete analysis of the few remaining cases.

\begin{prop}\label{prop-trans}
 For any set of matrices $\Sigma=\{A_0,A_1\} \subset \mathbb{R}^{2 \times 2},$ we have
 \begin{itemize}
 \item $\rho(\{A_0,A_1\})=\rho(\{A_0^{T},A_1^{T}\}),$ where $A^{T}$ is the transpose of $A$,
 \item $\rho(\{A_0,A_1\})=\rho(\{S A_0 S, S A_1 S\}),$ where $S= \begin{pmatrix} 0&1\\1&0 \end{pmatrix} $.
 \end{itemize}
 Moreover, in both cases the finiteness property holds for one set if and only if it holds for the other.\end{prop}

\begin{prop}\label{prop-sym}\cite[Proposition 6.13]{theys_thesis}
The finiteness property holds for sets of symmetric matrices.\end{prop}
\begin{proof}
The matrix norm induced by the euclidean vector norm is given by the largest singular value of the matrix. For symmetric matrices the largest singular value is equal to the largest magnitude of the eigenvalues. Thus $\max{\{||A||:A\in \Sigma\}}=\max{\{\rho(A):A\in \Sigma\}}$ and from  (\ref{a1}) it follows that $\rho(\Sigma)=\max{\{\rho(A):A \in \Sigma\}}$.
\end{proof}

\begin{prop}\label{prop-rholeq1} \label{cor-subid} \label{prop-A0leqA1}
Let $\Sigma = \{A_0,A_1\}\in \mathbb{N}^{n\times n}.$ The finiteness property holds in the following situations:
\begin{enumerate}
\item  $\rho (\Sigma)\leq 1,$
\item $A_0\leq I$ (or $A_1\leq I$).
\end{enumerate}
\end{prop}
\begin{proof} (1) It is known that for sets of nonnegative integer matrices, if $\rho\leq 1,$ then either $\rho=0$ and the finiteness property holds, or $\rho=1$, and there is a product of matrices in $\Sigma$ with a diagonal entry equal to one \cite{JungersProtasovBlondel06}.  Such a product $A\in \Sigma^t$ satisfies  $\rho(\Sigma)=\rho(A)^{1/t}=1$ and so the finiteness property holds when $\rho (\Sigma)\leq 1$.

(2) If $\rho(A_1) \leq 1$, then $\rho(A)\leq 1$ for all $A \in \Sigma^t$ and thus $\rho(\Sigma)\leq 1$ and the result follows from (1). If $\rho(A_1) >1$ then $\rho(\Sigma)=\rho(A_1)$ and so the finiteness property holds. 
\end{proof}

\begin{prop}\label{prop-switch}
Let $\Sigma = \{A_0,A_1\}\in \mathbb{N}^{n\times n}.$  The finiteness property holds in the following situations:
\begin{enumerate}
\item $A_0 \leq A_1$ (or $A_1 \leq A_0$),
\item $A_0 A_1\leq A_1^2$ (or $A_1 A_0\leq A_1^2$),
\item $A_0 A_1\leq A_1A_0$.
\end{enumerate}
\end{prop}

\begin{proof}

(1) Any product of length $t$ is bounded by  $A_1^t$.  Hence the joint spectral radius of $\Sigma$ is given by $\lim_{t\rightarrow \infty} ||A_1^t||^{1/t}=\rho(A_1)$.

(2) and (3). Let $A \in \Sigma^t$ be some product of length $t$. If $A_0 A_1\leq A_1^2$ or $A_0 A_1\leq A_1A_0$  we have $A \leq A_1^{t_1}A_0^{t_0}$ for some $t_0+t_1=t$.  The joint spectral radius is thus given by
\begin{eqnarray}\nonumber \rho = \lim_{t \rightarrow \infty}{\max_{t_1+t_0=t}{||A_1^{t_1}A_0^{t_0}||^{1/t}}}&\leq&\lim_{t \rightarrow \infty}{\max_{t_1+t_0=t}{||A_1^{t_1}||^{1/t}||A_0^{t_0}||^{1/t}}} \\ \nonumber &\leq&\max{(\rho(A_0), \rho(A_1))}.\end{eqnarray}
Hence the joint spectral radius is given by $\max{(\rho(A_0), \rho(A_1))}$.
\end{proof}

In order to analyse all possible sets of matrices, we consider all possible couples $(n_0,n_1)$, where $n_i$ is the number of nonzero entries in $A_i$.  From Proposition \ref{prop-switch}, we can suppose  $n_i=1,2,$ or $3$ and without loss of generality we take $n_0\leq n_1$.
\begin{itemize}
\item $n_0=1:$ \begin{itemize}
\item If $n_1=1$ or $n_1=2,$ the maximum row sum or the maximum column sum is equal to one for both matrices, and since these quantities are norms it follows from (\ref{a1}) that the joint spectral radius is less than one and from Proposition \ref{prop-rholeq1} that the finiteness property holds.
    \item If $n_1=3,$ it follows from Proposition \ref{prop-switch} that the only interesting cases are:  $$\Sigma=\left \{\begin{pmatrix}1&0\\0&0 \end{pmatrix},\begin{pmatrix}0&1\\1&1 \end{pmatrix}  \right \} \mbox{ and } \Sigma_0=\left \{\begin{pmatrix}0&1\\0&0 \end{pmatrix},\begin{pmatrix}1&0\\1&1 \end{pmatrix}  \right \}.$$ In the first case the matrices are symmetric and so the finiteness property holds by Proposition \ref{prop-sym}.  We keep $\Sigma_0$ for later.
\end{itemize}

\item $n_0=2:$ \begin{itemize}
\item $n_1=2:$ The only interesting cases are: $$\Sigma=\left \{\begin{pmatrix}1&1\\0&0 \end{pmatrix},\begin{pmatrix}0&1\\0&1 \end{pmatrix}  \right \} \mbox{ and } \Sigma_1=\left \{\begin{pmatrix}1&1\\0&0 \end{pmatrix},\begin{pmatrix}1&0\\1&0 \end{pmatrix}  \right \}.$$ Indeed in all the other cases either the maximum row sum or the maximum column sum is equal to one and the finiteness property follows from Proposition \ref{prop-rholeq1}.  The joint spectral radius of the first set is equal one.  Indeed, it is known that for nonnegative integer matrices, if the joint spectral radius is larger than one, then there must be a product of matrices with a diagonal entry larger than one  \cite{JungersProtasovBlondel06}.  This is impossible here, since as soon as a path leaves the first node, it cannot come back to it, and no path can leave the second node. By Proposition \ref{prop-rholeq1} the finiteness property holds for the first set. We keep $\Sigma_1$ for further analysis.
\item $n_1=3: $ If the zero entry of $A_1$ is on the diagonal (say, the second diagonal entry), then, from Proposition \ref{prop-A0leqA1}  we only need to consider the following case:
    $$\left \{\begin{pmatrix}0&1\\0&1 \end{pmatrix},\begin{pmatrix}1&1\\1&0 \end{pmatrix} \right \}.$$
    These matrices are such that $A_0A_1\leq A_1^2$ and so the finiteness property follows from Proposition \ref{prop-switch}.
    
    If the zero entry of $A_1$ is not a diagonal entry, we have to consider the following cases: $$ \Sigma_2=\left \{\begin{pmatrix}1&0\\1&0 \end{pmatrix} ,\begin{pmatrix}1&1\\0&1 \end{pmatrix} \right \} \mbox{ and } \Sigma_3=\left \{\begin{pmatrix}0&1\\1&0 \end{pmatrix}, \begin{pmatrix}1&1\\0&1 \end{pmatrix} \right \}.$$  We will handle $\Sigma_2$ and $\Sigma_3$ later.
\end{itemize}
\item $n_0,n_1=3:$  It has already been noticed by several authors (see, e.g., \cite[Proposition 5.17]{theys_thesis}) that $$\rho \left (\left \{\begin{pmatrix}1&1\\0&1 \end{pmatrix},\begin{pmatrix}1&0\\1&1 \end{pmatrix}  \right \}\right )=   \rho \left(\begin{pmatrix}1&1\\0&1 \end{pmatrix} \cdot  \begin{pmatrix}1&0\\1&1 \end{pmatrix} \right)^{1/2}  =\sqrt{\frac{1+\sqrt{5}}{2}}.$$  After excluding the case of symmetric matrices  and using the symmetry argument of Proposition \ref{prop-trans}, the only remaining case is:  $$ \left \{\begin{pmatrix}1&1\\0&1 \end{pmatrix},\begin{pmatrix}1&1\\1&0 \end{pmatrix}  \right \},$$ but again these matrices are such that $A_0A_1\leq A_1^2$ and so the finiteness property follows from Proposition \ref{prop-switch}.
\end{itemize}

We now analyse the cases that we have identified above.
For $\Sigma_0$, notice that $A_0^2\leq A_0 A_1$. Therefore, any product of length $t$ is dominated by a product of the form $A_1^{t_1}A_0A_1^{t_2}A_0\dots A_1^{t_l}$ for some  $t_1, t_l\geq 0$ and $t_i\geq 1$ $(i=2, \ldots, l-1)$. The norm of such a product is equal to $(t_1+1)(t_l+1)t_2\dots t_{l-1}$.  The maximal rate of growth of this norm with the product length is given by $\sqrt[5]{4}$   and so the joint spectral radius is equal to $\sqrt[5]{4}=\rho{(A_1^{4}A_0)}^{1/5}$. The maximal rate of growth is obtained for $t_i=4$.

For $\Sigma_1$,  simply notice that $\max_{A\in \Sigma^2}{\rho(A)}=\max_{A\in \Sigma^2}{||A||_\infty}=2,$ where $||\cdot||_\infty$ denotes the maximum row sum norm.  Hence by (\ref{a1}) we have $\rho(\Sigma)=\rho(A_0A_1)^{1/2}=\sqrt{2}$.

Consider now $\Sigma_2$. These matrices are such that $A_0^2 \leq A_0 A_1$  and so any product of length $t$ is dominated by a product of the form $A_1^{t_1}A_0A_1^{t_2}A_0\dots A_1^{t_l}$ for some  $t_1, t_l\geq 0$ and $t_i\geq 1$ $(i=2, \ldots, l-1)$. We have $$A_1^{t_1}A_0\dots A_1^{t_l}A_0=\begin{pmatrix}(t_1+1)\dots (t_l+1)&0\\(t_2+1)\dots (t_l+1)&0 \end{pmatrix}$$  and the maximum rate of growth of the norm of such a product is equal to $\sqrt{2}$. This rate is obtained for $t_i=3$  and $\rho=\rho(A_1^{3}A_0)^{1/4}=\sqrt{2}.$

The last case, $\Sigma_3$, is more complex and we give an independent proof for it.

\begin{prop}
The finiteness property holds for the set
$$\left \{\begin{pmatrix}1&1\\0&1 \end{pmatrix},\begin{pmatrix}0&1\\1&0 \end{pmatrix}  \right \}.$$
\end{prop}

\begin{proof}
Because $A_0^2=I$ we can assume the existence of a sequence of maximal-normed products $\Pi_i$ of length $L_i$, of the shape $B_{t_1} \dots B_{t_l}$ with $B_{t_i}=A_1^{t_i}A_0, \, \sum{t_k}+l=L_i,$ and $\lim{||\Pi_i||^{1/L_i}}=\rho(\Sigma)$.  We show that actually any maximal-normed product only has factors $B_3,$ except a bounded number of factors that are equal to $B_1,B_2,$ or $B_4$ and so the finiteness property holds.\\  Let us analyse one of these products $\Pi$. We suppose without loss of generality that $\Pi$ begins with a factor $B_3$.  First, it does not contain any factor $B_t:t>4$ because for such $t, \ B_{t-3}B_2\geq B_t$ and we can replace these factors without changing the length.\\  Now, our product $\Pi$ has less than $8$ factors $B_4$, because replacing the first seven factors $B_4$ with $B_3$, and the eighth one with $(B_3)^{3}$ we get a bigger-normed product of the same length (this is because $B_3\geq (3/4) B_4$, and $(B_3)^3\geq  (33/4)B_4$).  We remove these (at most) seven factors $B_4$ and by doing this, we just divide the norm by at most a constant $K_0$.\\ We now construct a bigger-normed product $\Pi'$ by replacing the left hand sides of the following inequalities by the respective right hand sides, which are products of the same length:
\begin{eqnarray} \nonumber B_i B_1B_1B_j\leq  B_i B_3B_j\\
\nonumber B_2 B_1B_2\leq  B_3 B_3\\
\nonumber B_3 B_1B_2\leq  B_2 B_2B_2\\
\nonumber B_2 B_1B_3\leq  B_2 B_2B_2\\
\nonumber \end{eqnarray}
If the factor $B_3B_1B_3$ appears eight times, we replace it seven times with $ B_2^3 \geq(4/5)B_3 B_1B_3$ and the last time with $B_2^3B_3^2$ which is greater than $7 B_2^3$.  By repeating this we get a new product $\Pi''\geq 7(4/5)^8\Pi'(1/K_0)>\Pi'(1/K_0)$ that has a bounded number of factors $B_1$.  We remove these factors from the product and by doing this we only divide by at most a constant $K_1$.\\
If there are more than four factors $B_2$ in the product, we replace the first three ones with $B_3$, and remove the fourth one.  It appears that for any $X\in \{B_2,B_3\},B_3^2X> 1.35B_3B_2X, $ and on the other hand, $B_3^2X\geq B_3^2B_2X\frac{1}{2,4349}$.  Then each time we replace four factors $B_2$ we get a new product: $\Pi'''\geq\frac{1.35^3}{2.4348}\Pi''(1/K_1)>\Pi''(1/K_1)$.  Finally we can remove the (at most) three last factors $B_2$ and by doing this, we only divide the product by at most a constant $K_2$.  By doing these operations to every $\Pi_i$, we get a sequence of products $\Pi_i'''$, of length at most $L_i$.  Now, introducing $K=K_0K_1K_2$, we compute $$\rho\geq \lim{||\Pi_i'''||^{1/(L_i)}}\geq \lim{||(1/K)\Pi_i||^{1/(L_i)}}=\rho.$$
Hence $\rho=\lim{||(A_1^3A_0)^t||^{1/(4t)}}=\rho(A_1^3A_0)^{1/4}=((3+\sqrt{13})/2)^{1/4}$, and the finiteness property holds.
\end{proof}
\end{section}
\begin{section}{Conclusion}
This paper provides a contribution to the analysis of the finiteness property for matrices that have rational entries. We have shown that the finiteness property holds for matrices with nonnegative rational entries if and only if it holds for pairs of matrices with binary entries. For pairs of binary matrices of dimension $2 \times 2$ we have shown that the property holds true and we conjecture that it holds for binary matrices of arbitrary dimension. A natural way to prove the conjecture for pairs of binary matrices would be to use induction on the size of the matrices, but this does not seem to be easy.  If the conjecture is true, it follows that the stability question for matrices with nonnegative rational entries is algorithmically decidable. If the conjecture is false, then the results and techniques developed in this paper can possibly help for constructing a counterexample.

\end{section}


\begin{thebibliography}{10}

\bibitem{berger-wang}
M.~A. Berger and Y.~Wang.
\newblock Bounded semigroups of matrices.
\newblock {\em Linear Algebra and its Applications}, 166:21--27, 1992.

\bibitem{blondel03undecidable}
V.~D. Blondel and V.~Canterini.
\newblock Undecidable problems for probabilistic automata of fixed dimension.
\newblock {\em Theory of Computing Systems}, 36(3):231--245, 2003.

\bibitem{BlondelJungersProtasov06}
V.~D. Blondel, R.~Jungers, and V.~Protasov.
\newblock On the complexity of computing the capacity of codes that avoid
  forbidden difference patterns.
\newblock {\em IEEE Transactions on Information Theory}, 52(11):5122--5127,
  2006.

\bibitem{blondel-elementary}
V.~D. Blondel, J.~Theys, and A.~Vladimirov.
\newblock An elementary counterexample to the finiteness conjecture.
\newblock {\em SIAM Journal on Matrix Analysis}, 24(4):963--970, 2003.

\bibitem{tsitsiklis97lyapunov}
V.~D. Blondel and J.~N. Tsitsiklis.
\newblock The lyapunov exponent and joint spectral radius of pairs of matrices
  are hard - when not impossible - to compute and to approximate.
\newblock {\em Mathematics of Control, Signals, and Systems}, 10:31--40, 1997.
\newblock Correction in 10, 381.

\bibitem{blondel00boundedness}
V.~D. Blondel and J.~N. Tsitsiklis.
\newblock The boundedness of all products of a pair of matrices is undecidable.
\newblock {\em Systems and Control Letters}, 41(2):135--140, 2000.

\bibitem{cfbousch}
T.~Bousch and J.~Mairesse.
\newblock Asymptotic height optimization for topical ifs, tetris heaps, and the
  finiteness conjecture.
\newblock {\em Journal of the Mathematical American Society}, 15(1):77--111,
  2002.

\bibitem{Crespi05}
V.~Crespi, G.~V. Cybenko, and G.~Jiang.
\newblock {The Theory of Trackability with Applications to Sensor Networks}.
\newblock Technical Report TR2005-555, Dartmouth College, Computer Science,
  Hanover, NH, August 2005.

\bibitem{Gu1}
L.~Gurvits.
\newblock Stability of discrete linear inclusions.
\newblock {\em Linear Algebra and its Applications}, 231:47--85, 1995.

\bibitem{master_thesis}
R.~Jungers.
\newblock {NP}-completeness of the computation of the capacity of constraints
  on codes.
\newblock Master's thesis, Universit\'e Catholique de Louvain, 2005.

\bibitem{JungersProtasovBlondel06}
R.~Jungers, V.~Protasov, and V.~D. Blondel.
\newblock Efficient algorithms for deciding the type of growth of products of
  integer matrices.
\newblock submitted to publication.

\bibitem{cfkoz}
V.~Kozyakin.
\newblock A dynamical systems construction of a counterexample to the
  finiteness conjecture.
\newblock {\em Proc. 44th IEEE Conference on Decision and Control and ECC
  2005}, 2005.

\bibitem{lagarias-finiteness}
J.~C. Lagarias and Y.~Wang.
\newblock The finiteness conjecture for the generalized spectral radius of a
  set of matrices.
\newblock {\em Linear Algebra and its Applications}, 214:17--42, 1995.

\bibitem{lib}
D.~Liberzon, J.~P. Hespanha, and A.~S. Morse.
\newblock Stability of switched systems: a lie-algebraic condition.
\newblock {\em Systems and Control Letters}, 37(3):117--122, 1999.

\bibitem{moision-joint-constraints}
B.~E. Moision, A.~Orlitsky, and P.~H. Siegel.
\newblock On codes with local joint constraints.
\newblock Unpublished.

\bibitem{mos-bounds}
B.~E. Moision, A.~Orlitsky, and P.~H. Siegel.
\newblock Bounds on the rate of codes which forbid specified difference
  sequences.
\newblock In {\em Proc. 1999 IEEE Global Telecommun. Conf. (GLOBECOM'99)},
  1999.

\bibitem{MoisionOrlitskySiegel01}
B.~E. Moision, A.~Orlitsky, and P.~H. Siegel.
\newblock On codes that avoid specified differences.
\newblock {\em IEEE Transactions on Information Theory}, 47:433--442, 2001.

\bibitem{protasov3}
V.~Y. Protasov.
\newblock On the asymptotics of the partition function.
\newblock {\em Sb. Math.}, 191(3-4):381--414, 2000.

\bibitem{theys_thesis}
J.~Theys.
\newblock {\em Joint Spectral Radius : theory and approximations}.
\newblock PhD thesis, Universit\'e Catholique de Louvain, 2005.

\end{thebibliography}
\end{document}